\documentclass[preprint,12pt]{elsarticle}

\usepackage{amssymb}

\journal{arXiv}

\begin{document}

\begin{frontmatter}

\title{Multivariate two-sample extended empirical likelihood}
\author{Fan Wu\footnote{Corresponding author; email address: fwu@uvic.ca} and Min Tsao}
\address{Department of Mathematics and Statistics, University of Victoria, Victoria, British Columbia, Canada V8W 3R4}

\begin{abstract}

Jing (1995) and Liu et al. (2008) studied the two-sample empirical likelihood and showed it is Bartlett correctable for the univariate and multivariate cases, respectively. We expand its domain to the full parameter space and obtain a two-sample extended empirical likelihood which is more accurate and can also achieve the second-order accuracy of the Bartlett correction.
\\[0.2in]
\end{abstract}


\end{frontmatter}

\section{Introduction}  
The empirical likelihood introduced by Owen (1988, 1990) is a versatile non-parametric method of inference with many applications (Owen, 2001). One problem which the empirical likelihood method has been successfully applied to is the two-sample problem (Jing, 1995; Liu et al. 2008; Wu and Yan, 2012) where the parameter of interest  $\theta$ is the difference between the means of two populations. The well-known Behrens-Fisher problem is a special two-sample problem where the two populations are known to be normally distributed.  Following DiCiccio et al. (1991) who showed the surprising result that the (one-sample) empirical likelihood for a smooth function of the mean is Bartlett correctable, Jing (1995) and Liu et al. (2008) proved that the two-sample empirical likelihood for $\theta$ is also Bartlett correctable for the univariate and multivariate cases, respectively. The coverage error of a confidence region based on the original empirical likelihood is $O(n^{-1})$, but that based on the Bartlett corrected empirical likelihood is only $O(n^{-2})$.

For a one-sample empirical likelihood, there is a mismatch between its domain and the parameter space in that it is defined on only a part of the parameter space. This mismatch is a main cause of the undercoverage problem associated with empirical likelihood confidence regions (Tsao, 2013). The two-sample empirical likelihood for $\theta$ also has the mismatch problem as it is defined on a bounded region but the parameter space is $\mathbb{R}^{d}$. In this paper, we derive an extended version of the original two-sample empirical likelihood (OEL) by expanding its domain into $\mathbb{R}^{d}$ through the composite similarity mapping of Tsao and Wu (2013). The resulting two-sample extended  empirical likelihood (EEL) for $\theta$ is defined on the entire $\mathbb{R}^{d}$ and hence free from the mismatch problem. Under mild conditions, this EEL has the same asymptotic properties as the OEL. It can also attain the second order accuracy of the two-sample Bartlett corrected empirical likelihood (BEL)  of Jing (1995) and Liu et al. (2008). The first order version of this EEL is substantially more accurate than the OEL, especially for small sample sizes. It is also easy to compute and competitive in accuracy to the second order BEL. We recommend it for two-sample empirical likelihood inference.

\section{Two-sample empirical likelihood}

Let $\{X_1, \dots, X_m\}$ and $\{Y_1, \dots, Y_n\}$ be independent copies of random vectors $X\in \mathbb{R}^{d}$ and $Y\in \mathbb{R}^{d}$, respectively. Denote by $\mu_x$ and $\Sigma_x$ the mean and covariance matrix of $X$, and by $\mu_y$ and $\Sigma_y$ the mean and covariance matrix of $Y$, respectively. The unknown parameter of interest is the difference in means $\theta_0=\mu_y-\mu_x\in \mathbb{R}^{d}$ and the parameter space is the entire $\mathbb{R}^{d}$. We will need the following three conditions later in the paper:\vspace{0.1in}

\noindent{\em C1.} $\Sigma_x$ and $\Sigma_y$ are finite covariance matrix with full rank $d$;

\noindent{\em C2.} $\lim \sup_{\|t\|\rightarrow \infty}|E[\exp\{it^TX\}]| <1$ and $\lim \sup_{\|t\|\rightarrow \infty}|E[\exp\{it^TY\}]| <1$;

\noindent{\em C3.} $E\|X\|^{15}<+\infty$ and  $E\|Y\|^{15}<+\infty$. \\[0.1in]
Denote by $p=(p_1,...,p_m)$ and $q=(q_1,...,q_n)$ two probability vectors satisfying $p_i\geq 0$, $q_j\geq 0$, $\sum_{i=1}^m p_i=1$ and $\sum_{i=1}^n q_j=1$. Let $\mu_x(p)=\sum_{i=1}^m p_i X_i$ and  $\mu_y(q)=\sum_{j=1}^n q_j Y_j$, and denote by $\theta(p,q)$ their difference, that is, $$\theta(p,q)=\mu_y(q)-\mu_x(p).$$
The original two-sample empirical likelihood for a $\theta \in \mathbb{R}^d$, $L(\theta)$, is defined as
\begin{equation}
L(\theta) = \max_{(p,q):\theta(p,q)=\theta} \left(\prod^{m}_{i=1}p_i \right) \left(\prod^{n}_{j=1}q_i \right). \label{oel}
\end{equation}
The corresponding two-sample empirical log-likelihood ratio for $\theta$ is thus
\begin{equation}
l(\theta) = -2\max_{(p,q):\theta(p,q)=\theta} \left(\sum^{m}_{i=1}\log(mp_i)+ \sum^{n}_{j=1}\log(nq_i) \right). \label{oel.log}
\end{equation}

In order to develop our extended empirical likelihood, it is important to understand the domains of the original empirical likelihood ratio $L(\theta)$ and log-likelihood ratio $l(\theta)$. The domain of $L(\theta)$ is given by
\begin{center}
$D_\theta=\{\theta\in \mathbb{R}^d: \mbox{there exist $p$ and $q$ such that $\mu_x(p)=\sum_{i=1}^m p_i X_i$,}$ \\ 
                  $\mbox{$\mu_y(q)=\sum_{j=1}^n q_j Y_j$ and $\theta=\theta(p,q)=\mu_y(q)-\mu_x(p)$}\}.$
\end{center}
Since the range of $\mu_x(p)$ and $\mu_y(q)$ are the convex hulls of the $X_i$ and $Y_i$, respectively, $D_\theta$ is a bounded, closed and connected region in $\mathbb{R}^d$ without voids. More details about the geometric structure of $D_\theta$ may be found in Lemma 1 and its proof, and one of these is that an interior point of $D_\theta$ can be expressed as $\theta(p,q)=\mu_y(q)-\mu_x(p)$ where all elements of $p$ and $q$ are straightly positive. Correspondingly, a boundary point of $D_\theta$ can only be expressed as $\theta(p,q)=\mu_y(q)-\mu_x(p)$ where one or more elements of $p$ and $q$ are zero. This implies that $L(\theta)=0$ for any boundary point $\theta$ of $D_\theta$ and $L(\theta)>0$ for every interior point. We define the domain of the empirical log-likelihood ratio $l(\theta)$ as
\[ 
\Theta_n=\{\theta: \theta \in D_\theta \hspace{0.1in} \mbox{and} \hspace{0.1in} l(\theta)<+\infty \},
\]
which excludes the boundary points of $D_\theta$. That is, $\Theta_n$ is the open set in $\mathbb{R}^d$ that contains the collection of $\theta(p,q)=\mu_y(q)-\mu_x(p)$ where elements of $p$ and $q$ are all straightly positive.
To differentiate between the $l(\theta)$ in (\ref{oel.log}) and the extended version of $l(\theta)$ in the next section, we will refer to the $l(\theta)$ in (\ref{oel.log}) as the original two-sample empirical log-likelihood ratio or simply ``OEL $l(\theta)$''. The extended version will be referred to as the ``EEL $l^*(\theta)$''. 

Let $N=m+n$, $f_m=N/m$ and $f_n=N/n$. Without loss of generality, assume that $m\geq n>d$.
By the method of Lagrangian multipliers, we have
\begin{equation}
l(\theta) = 2\left[\sum^{m}_{i=1}\log \{1-f_m \lambda^{T}(X_i-\mu_x)\}+ \sum^{n}_{j=1} \log \{1+f_n \lambda^{T}(Y_j-\mu_y)\} \right] \label{oel.lag.mul}
\end{equation}
where the multiplier $\lambda=\lambda(\theta)$ satisfies
\begin{equation}
   \sum^{m}_{i=1} \frac{X_i-\mu_x}{1-f_m \lambda^{T}(X_i-\mu_x)}=0 \hspace{0.2in} \mbox{and} \hspace{0.2in}
    \sum^{n}_{j=1} \frac{Y_j-\mu_y}{1+f_n \lambda^{T}(Y_j-\mu_y)}=0,
\end{equation}
and
\begin{equation}
   \sum^{n}_{j=1} \frac{Y_j}{1+f_n \lambda^{T}(Y_j-\mu_y)}-\sum^{m}_{i=1} \frac{X_i}{1-f_m \lambda^{T}(X_i-\mu_x)}=\theta.    
\end{equation}
Under the assumption ($C1$), Jing (1995) and Liu et al. (2008) showed that 
\begin{equation}
l(\theta_0) \stackrel{D}{\longrightarrow}\chi^2_d \hspace{0.2in} \mbox{as} \hspace{0.1in} n\rightarrow +\infty.  \label{oel.chi}
\end{equation}
Hence, the 100($1-\alpha$)\% OEL confidence interval for $\theta_0$ is
\begin{equation}
{\mathcal C}_{1-\alpha}=\{\theta: \theta \in \mathbb{R}^{d} \mbox{ and } l(\theta) \leq c_\alpha \}   \label{oel.interval}
\end{equation}
where $c_\alpha$ is ($1-\alpha$)th quantile of the $\chi^2_d$ distribution. The coverage error of $\mathcal{C}_{1-\alpha}$ is $O(n^{-1})$, that is
\begin{equation}
P(\theta_0\in {\mathcal C}_{1-\alpha})= P(l(\theta_0) \leq c_\alpha)=1-\alpha + O(n^{-1}). \label{oel.coverage}
\end{equation}
Under assumptions ($C1$), ($C2$) and ($C3$), Jing (1995) and Liu et al. (2008) also showed that the OEL $l(\theta)$ is Bartlett correctable, that is
\begin{equation}
P(\theta_0\in {\mathcal C}_{1-\alpha}')=P(l(\theta_0)\leq c_\alpha (1+\eta N^{-1}))= \alpha + O(n^{-2}) \label{oel.bart}
\end{equation}
where ${\mathcal C}_{1-\alpha}'= \{\theta: l(\theta) \leq c_\alpha (1+\eta N^{-1}) \}$ is the Bartlett corrected empirical likelihood (BEL) confidence interval and $\eta$ is the Bartlett correction constant in Theorem 2 of Liu et al. (2008). For the one-dimensional case, a formula for this constant was first given in Theorem 2 in Jing (1995) but the formula is incomplete (Liu et al. 2008).  
See also Wu and Yan (2012) and Qin (1994) for discussions about two-sample empirical likelihood methods.

\section{Two-sample extended empirical likelihood}

Like the one-sample empirical likelihood for the mean, the two-sample OEL $l(\theta)$ also suffers from the mismatch problem between its domain and the parameter space since the parameter space is $\mathbb{R}{^d}$ but $\Theta_n\subset \mathbb{R}{^d}$. This is a main cause of the undercoverage problem of empirical likelihood confidence regions (Tsao, 2013; Tsao and Wu, 2013). 
To overcome the mismatch, we now extend the OEL $l(\theta)$ by expanding its domain to the entire $\mathbb{R}{^d}$. 

For simplicity, in addition to $m\geq n$ we further assume that $m/n=O(1)$ so that $O(n^{-1})$, $O(m^{-1})$ and $O(N^{-1})$, for example, are all interchangeable. A point estimator for $\theta_0$ is $\hat{\theta}=\bar{Y}-\bar{X}$ where $\bar{X}=m^{-1} \sum X_i$ and $\bar{Y}=n^{-1} \sum Y_j$ are the sample means. It is easy to verify that $\hat{\theta}$ is the maximum empirical likelihood estimator (MELE) for $\theta_0$. 
Following Tsao and Wu (2013), we define the composite similarity mapping $h_N^C: \Theta_n \rightarrow \mathbb{R}^{d}$ centred on $\hat{\theta}$ as
\begin{equation}
h_N^C(\theta)=\hat{\theta}+\gamma(N,l(\theta))(\theta-\hat{\theta}) \label{h.function}
\end{equation}
where function $\gamma(n,l(\theta))$ is the expansion factor given by 
\begin{equation}
\gamma(N,l(\theta))=1+\frac{l(\theta)}{2N}. \label{eel1.00}
\end{equation} 

To investigate the properties of the composite similarity mapping $h_N^C$, we need Lemma 1 below which gives two properties of the two-sample OEL $l(\theta)$. For convenience, we denote by $[\hat{\theta},\theta]$ the line segment that connects $\hat{\theta}$ and $\theta$ and by $\theta_b$ a boundary point of $\Theta_n$. We have
\vspace{0.2in}

\noindent \textbf{Lemma 1.} \emph{The two-sample OEL $l(\theta)$ satisfies: ($i$) if $\theta\in \Theta_n$ and $\theta' \in [\hat{\theta},\theta]$, then $l(\theta')\leq l(\theta)$ and ($ii$) for $\theta \in \Theta_n$, $\lim_{\theta \rightarrow \theta_b}l(\theta)=+\infty$.}
\vspace{0.2in}

\noindent Lemma 1 shows the two-sample OEL $l(\theta)$ for the difference of two means behaves exactly like its one-sample counterpart for the mean in terms monotonicity and boundary behaviour: it is ``monotone increasing'' along each ray originating from the MELE and it goes to infinity as $\theta$ approaches a boundary point from within $\Theta_n$. Nevertheless, the two-sample and  one-sample cases are not entirely the same; the contours of the two-sample OEL may not be convex but that of the one-sample OEL always are. Theorem 1 below gives three key properties of composite similarity mapping $h_N^C: \Theta_n \rightarrow \mathbb{R}^{d}$.\\[0.1in]

\noindent {\bf Theorem 1.} \emph{
Under the assumption (C1), $h^C_N: \Theta_n\rightarrow \mathbb{R}^{d}$ defined by (\ref{h.function}) and (\ref{eel1.00}) satisfies (i) it has a unique fixed point at $\hat{\theta}$, (ii) it is a similarity mapping for each individual contour of the OEL $l(\theta)$  and (iii) it is a bijective mapping from $\Theta_n$ to $\mathbb{R}^{d}$. }
\vspace{0.2in}

Since $h_N^{C}: \Theta_n\rightarrow \mathbb{R}^{d}$ is bijective, it has an inverse function which we denote by $h_N^{-C}(\theta): \mathbb{R}^{d} \rightarrow \Theta_n$. For any $\theta\in \mathbb{R}^{d}$, let $\theta'=h_N^{-C}(\theta) \in \Theta_n$.  The two-sample extended empirical log-likelihood ratio EEL $l^*(\theta)$ is given by
\begin{equation}
l^*(\theta)=l(h_N^{-C}(\theta))=l(\theta'), \label{eel1}
\end{equation}
which is defined for $\theta$ values throughout $\mathbb{R}^{d}$. Hence the EEL $l^*(\theta)$ is free from the mismatch problem of the OEL $l(\theta)$. Denote by $\theta_0'$ the image of $\theta_0$ under the inverse transformation $h_N^{-C}(\theta): \mathbb{R}^{d} \rightarrow \Theta_n$, that is 
\begin{equation}
h_N^{-C}(\theta_0)= \theta_0'. \label{temp3}
\end{equation}
Then, the EEL $l^*(\theta)$ evaluated at $\theta_0$ is given by
\begin{equation}
l^*(\theta_0) = l(h_N^{-C}(\theta_0)) = l(\theta'_0)=l(\theta_0+\theta'_0-\theta_0). \label{eel1.thm.01}
\end{equation}

If $|\theta'_0-\theta_0|$ is very small, then $l^*(\theta_0)$ will have the same asymptotic distribution as $l^*(\theta)$. Lemma 2 below shows that $\theta_0' \in [\hat{\theta},\theta_0]$ and that $|\theta'_0-\theta_0|$ is indeed very small. \\[0.1in]

\noindent \textbf{Lemma 2.}
{\em Under assumption (C1), point $\theta'_0$ defined by equation (\ref{temp3}) satisfies }
$(i)$ $\theta'_0 \in [\hat{\theta},\theta_0]$ and $(ii)$ $\theta'_0-\theta_0=O_p(n^{-3/2})$.

\vspace{0.2in}

Theorem 2 below shows that EEL $l^*(\theta_0)$ has the same asymptotic chi-square distribution as the OEL $l(\theta_0)$.

\vspace{0.2in}

\noindent \textbf{Theorem 2.} \emph{
Under assumption (C1), the EEL $l^*(\theta_0)$ defined by (\ref{eel1.thm.01}) satisfies
\begin{equation}
l^*(\theta_0) \stackrel{D}{\longrightarrow} \chi^2_d  \hspace{0.2in} as \hspace{0.1in} n\rightarrow +\infty. \label{eel1.thm}
\end{equation}
} \\
By Theorem 2, the $100(1-\alpha)\%$ EEL confidence interval for $\theta_0$ is 
\begin{equation}
{\mathcal C}^*_{1-\alpha}=\{\theta: \theta \in \mathbb{R}^{d} \mbox{ and } l^*(\theta) \leq c_\alpha \}, \label{eel.ci}
\end{equation}
which has a coverage error of $O(n^{-1})$. The expansion factor in (\ref{eel1.00}) is a convenient choice which also gives good numerical results. There are many other choices available under which Theorems 1 and 2 also hold. This provides an opportunity to optimize the choice of expansion factor to obtain the second order accuracy. Theorem 3 below gives such an optimal choice.

\vspace{0.2in}

\noindent \textbf{Theorem 3.} \emph{
Under assumptions (C1), (C2) and (C3), and let $l_2^*(\theta)$ be the EEL defined by the composite similarity mapping 
(\ref{h.function}) with the following expansion factor 
\begin{equation}
\gamma_2(N,l(\theta))=1+\frac{\eta}{2N}\left[l(\theta)\right]^{\delta(N)} \label{eel2.01}
\end{equation}
where $\delta(N)=O(n^{-1/2})$ and $\eta$ is the Bartlett correction factor for the two-sample empirical likelihood in (\ref{oel.bart}). Then, we have
\begin{equation}
l^*_2(\theta_0)= l(\theta_0)\left[1-\eta/N + O_p(n^{-3/2})\right], \label{eel2.thm}
\end{equation}
and
\begin{equation}
P(l^*_2(\theta_0)\leq c)= P(\chi^2_d \leq c) +O(n^{-2}). \label{eel2.thm.2}
\end{equation}
}

Replacing EEL $l^*(\theta)$ in (\ref{eel.ci}) with $l_2^*(\theta)$ gives an EEL confidence interval which, by (\ref{eel2.thm.2}), has a coverage error of $O(n^{-2})$. Because of this, we call $l_2^*(\theta)$ the {\em second order} EEL or EEL$_2$. Correspondingly, we call the EEL $l^*(\theta)$ defined by expansion factor (\ref{eel1.00}) the {\em first order} EEL or EEL$_1$.

To limit the length of this paper, we have not included the proofs of the above lemmas and theorems here. They are given in a technical report (Wu and Tsao, 2013) available on request from the authors.

\section{Numerical examples}

We now compare the coverage accuracy of 95\% confidence intervals based on the OEL, BEL and EEL with two numerical examples. Comparisons based on 90\% and 99\% confidence intervals give similar conclusions and are thus not included. They can also be found in Wu and Tsao (2013). In the following, $BVN(0,I)$ represents the standard bivariate normal distribution and $X\sim (\chi^2_1,\chi^2_1)^T$, for example, represents a bivariate random vector $X$ whose two elements are independent $\chi^2_1$ random variables.

\vspace{0.1in}
\noindent {\bf Example 1:} $X\sim (\chi^2_1, \chi^2_1)^T$ and $Y\sim BVN(0,I)$. \\
\noindent {\bf Example 2:} $X\sim (\chi^2_3, \chi^2_3)^T$ and $Y\sim (Exp(1), Exp(1))^T$. \\[0.1in]
To see the effect of the composite similarity mapping, Figure 1 compares contours for the OEL $l(\theta)$ and the corresponding contours for the EEL$_1$ $l(\theta)$ based on the same pair of $X$ and $Y$ samples from Example 1. We see that the contours in the two plots are identical in shape and the contours in both plots are centred on the MELE $\hat{\theta}$ as indicated in Theorem 1. Further, at any fixed level, the contour of the EEL $l^*(\theta)$ is larger in scale. 

Simulated coverage probabilities for the two examples are given in Tables 1 and 2, respectively. Each simulated probability in the tables is based on 10,000 pairs of random samples whose sizes are indicated by the row and column headings, respectively. The BEL and EEL$_2$ were computed by using the estimated Bartlett correction factor from page 550 in Liu et al. (2008).  We summarize the tables with the following observations: (1) EEL$_1$ is consistently more accurate than the OEL. Surprisingly, it is also more accurate than the second order BEL and EEL$_2$ for small and moderate sample sizes ($n,m \leq 20$) and competitive in accuracy when sample sizes are larger. (2) EEL$_2$ is more accurate than OEL and BEL for small and moderate sample sizes. It is comparable to BEL when one or both sample sizes are large. 

To conclude, EEL$_1$ is easy-to-compute and is the most accurate overall. Hence, we recommend EEL$_1$ for two-sample problems.

\begin{figure}
\includegraphics[height=3in,width=5in,angle=0]{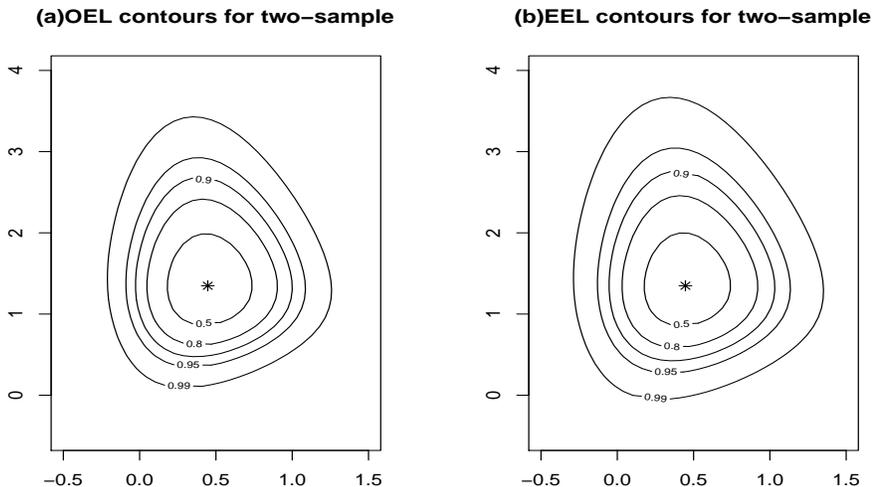}
\caption{(a) Two-sample OEL contours; (b) Two-sample EEL contours. Both plots are based the same pair of $X$ and $Y$ samples from Example 1 with sample size $n=20$ and $m=20$. The star in the middle of the plot is the MELE $\hat{\theta}$. EEL$_1$ contours are larger than but similar to OEL contours with the same centre and identical shape. }
\label{fig1}
\end{figure}

\begin{table}[htpb]
\begin{center}
{\bf \footnotesize Table 1: Coverage probabilities of 95\% OEL, EEL$_1$, BEL \& EEL$_2$ \\confidence intervals: $X\sim (\chi^2_1,\chi^2_1)$ and $Y\sim BVN(0,I)$} \\[0.1in]
{\footnotesize
\addtolength{\tabcolsep}{+5pt}
\begin{tabular}{c l c c c c  } \hline
&  &    $n$=10& $n$=20&$n$=30&$n$=40 \\  \hline
$m$=10  	          
          	&OEL        	&84.0	&88.7	&89.4	&89.8	\\	
          	&EEL$_1$     	&90.2	&92.4	&92.2	&92.0	\\	
          	&BEL        	&86.8	&90.7	&91.4	&91.6	\\	\vspace{0.07in}
          	&EEL$_2$    	&88.5	&91.7	&92.1	&92.1	\\	

$m$=20    
          	&OEL        	&83.0	&90.1	&91.8	&92.2	\\	
          	&EEL$_1$     	&87.3	&93.0	&93.8	&93.8	\\	
          	&BEL        	&85.3	&91.7	&93.1	&93.4	\\	\vspace{0.07in}
          	&EEL$_2$    	&86.4	&92.2	&93.5	&93.6	\\

$m$=30    
          	&OEL        	&80.2	&90.1	&91.9	&92.8	\\	
          	&EEL$_1$     	&83.7	&92.3	&93.6	&94.2	\\	
          	&BEL        	&82.7	&91.6	&93.0	&93.7	\\	\vspace{0.07in}
          	&EEL$_2$    	&83.7	&92.0	&93.3	&93.9	\\	
$m$=40   
          	&OEL        	&78.9	&88.3	&91.3	&92.4	\\	
          	&EEL$_1$     	&81.8	&90.2	&92.6	&93.7	\\	
          	&BEL        	&81.6	&89.8	&92.2	&93.4	\\	
          	&EEL$_2$    	&82.5	&90.3	&92.4	&93.6	\\	   	
                                                
\end{tabular}
}
\end{center} 
\end{table}

\begin{table}[htpb]
\begin{center}
{\bf \footnotesize Table 2: Coverage probabilities of 95\% OEL, EEL$_1$, BEL \& EEL$_2$ \\confidence intervals: $X\sim (\chi^2_3,\chi^2_3)$ and $Y\sim (Exp(1), Exp(1))$} \\[0.1in]
{\footnotesize
\addtolength{\tabcolsep}{+5pt}
\begin{tabular}{c l c c c c  } \hline
&  &    $n$=10& $n$=20&$n$=30&$n$=40 \\  \hline
$m$=10  	          
          	&OEL        	&81.8	&89.7	&91.0	&91.3	\\	
          	&EEL$_1$     	&88.8	&93.4	&93.8	&93.5	\\	
          	&BEL        	&84.6	&91.5	&92.6	&92.7	\\	\vspace{0.07in}
          	&EEL$_2$    	&87.0	&92.4	&93.2	&93.1	\\

$m$=20    
          	&OEL        	&81.2	&89.9	&91.8	&92.4	\\	
          	&EEL$_1$     	&85.6	&92.4	&93.8	&94.0	\\	
          	&BEL        	&83.7	&91.4	&92.9	&93.5	\\	\vspace{0.07in}
          	&EEL$_2$    	&85.2	&91.9	&93.2	&93.8	\\	
       	 
$m$=30    
          	&OEL        	&78.2	&89.1	&91.3	&92.9	\\	
          	&EEL$_1$     	&84.4	&91.6	&93.0	&94.4	\\	
          	&BEL        	&80.7	&90.7	&92.4	&93.8	\\	\vspace{0.07in}
          	&EEL$_2$    	&84.4	&91.4	&92.7	&94.0	\\

$m$=40   
          	&OEL        	&79.6	&89.5	&91.5	&92.5	\\	
          	&EEL$_1$     	&83.0	&91.5	&93.0	&93.6	\\	
          	&BEL        	&82.2	&91.1	&92.6	&93.3	\\	
          	&EEL$_2$    	&83.9	&91.6	&92.9	&93.5	\\	
                                       
\end{tabular}
}
\end{center} 
\end{table}

\newpage
\noindent {\bf References}






\end{document}